\documentclass[12pt, amsfonts, epsfig]{amsart}
\usepackage{latexsym, amssymb, amsmath, epsfig, amscd, amsfonts}
 \newlength{\baseunit}               % the basic unit length
 \newcount{\numlines}                % depth of picture (in number of lines)
\newcommand{\N}{\mathbb{N}}

\newcommand{\R}{\mathbb{R}}
\newcommand{\C}{\mathbb{C}}

\newcommand{\proj}{\mathbb P}

        \newfont{\hollow}{msbm10 scaled\magstep1}
        
        \newfont{\Bfmit}{eufm10 scaled\magstep1}

\newcommand{\PGL}{\operatorname{PGL}}

\newcommand{\Gr}{\operatorname{Gr}}

\newtheorem{thm}{Theorem}[subsection]

\newtheorem{cor}[thm]{Corollary}

\theoremstyle{definition}

\newtheorem{rem}[thm]{Remark}           %\renewcommand{\theremark}{}
\newtheorem{toric}[thm]{Toric Setting}            %\renewcommand{\thenote}{}

\theoremstyle{remark}

\newcommand{\lremind}[1]{{}}
\newcommand{\bremind}[1]{{}}

\newcommand{\cut}[1]{}
\begin{document}
\pagestyle{plain} \title{{ \large{Toric Degenerations of GIT
Quotients, Chow Quotients, and $\overline{M}_{0,n}$} } }
\author{Yi Hu }

\address{Department of Mathematics, University of Arizona, Tucson,
AZ 86721}

\email{yhu@math.arizona.edu}

\maketitle

%\tableofcontents

\section{Introduction }

The moduli space $\overline{M}_{0,n}$ plays important roles in
algebraic geometry and theoretical physics. Yet, some basic
properties of $\overline{M}_{0,n}$ still remain open. For example,
$\overline{M}_{0,n}$ is rational and nearly toric (that is, it
contains  a toric variety as a Zariski open subset), but it is not
a toric variety itself starting from dimension 2 ($n \ge 5$). So,
a basic question is: {\it Can it be degenerated flatly to a
projective toric variety}?  Finding toric degeneration is
important because many calculations may then be done in terms of
combinatorial data extracted out of  polytopes and/or fans. The
main purpose of this note is to answer the above question in
affirmative.

We achieve this  by showing that the Chow quotient of the
Grassmannian $\Gr(2, {\mathbb C}^n)$ admits toric degeneration,
which in turn, follows from a  theorem that we prove for toric
degenerations of more general Chow quotients. Consequently, this
implies as well that the Chow quotients of higher Grassmannians
(\cite{Ka}, \cite{Lafforgue99}, \cite{Lafforgue}), hence also (the
main component of) the moduli spaces of hyperplane arrangements
(\cite{HKT}) and (the main component of) the moduli spaces of semi
log canonical pairs (\cite{Alexeev96}, \cite{Alexeev02},
\cite{Alexeev06}, \cite{Hacking03}, \cite{KT}), can be degenerated
flatly to projective toric varieties. Along the way, we also
argued that GIT quotient of a flat family is again flat. In
particular, all GIT quotients of flag varieties by maximal tori
can be flatly degenerated  to  projective toric varieties
(\cite{FH} and \cite{HMSV05}).

\section{Flat family of GIT quotients}

\subsection{ Flat family of $G$-varieties.}
Throughout the paper, we will work over the field of complex
numbers.  Let $f :Y \rightarrow B$ be a flat family of projective
varieties. Assume that a reductive algebraic group $G$ acts
algebraically on the family and preserves fibers. That is,  $G$
acts on the total space $Y$, acts on the base trivially, and $f$
is $G$-equivariant.

Let $L$ be any $G$-linearized ample line bundle over $Y$ and $L_t$
be the restriction of $L$ to $Y_t= f^{-1}(t)$ with the induced
linearization for any $t \in B$.  By \cite{GIT}, we have
$$Y_t^{ss}(L_t) = Y^{ss}(L) \cap Y_t$$ and
$$Y_t^s(L_t) = Y^s(L) \cap Y_t.$$
Using these identities, the assumption that $G$ acts trivially on
$B$, and $f$ is $G$-equivariant, we obtain an induced map
$$\underline{f} :Y^{ss}(L)/\!/G \rightarrow B$$
with fibers $Y_t^{ss}(L_t)/\!/G$ for any $t \in B$.

When $B$ is regular and of dimension 1, the flatness of
$\underline{f}$ is automatic. We will show that under the above
conditions, it is always flat.

\begin{thm}
\label{flatgit} The induced family $\underline{f} :Y^{ss}(L)/\!/G
\rightarrow B$ of GIT quotients $Y_t^{ss}(L_t)/\!/G$ ($t \in B$)
is always flat.
\end{thm}

\proof Replace $L$ by a large tensor power, we may assume that $L$
is very ample and descends to a very ample line bundle $M$ over
the quotient $Y^{ss}(L)/\!/G$. Hence for every $n >0$,  $L^n$
descends to the line bundle $M^n$ over $Y^{ss}(L)/\!/G$. It
follows that $L_t^n$ descends to $M_t^n$ for all $t \in B$
simultaneously. Now, by the Theorem of  $``$Quantization commutes
with Reduction$"$ (see, for example, \cite{Teleman}), we have that
for $n >>0$,
$$H^0(Y_t^{ss}(L_t)/\!/G, M_t^n) = H^0(Y_t, L_t^n)^G.$$
This shows that $\dim H^0(Y_t^{ss}(L_t)/\!/G, M_t^n)$ is
independent of the parameter $t$ because  $\dim H^0(Y_t, L_t^n)^G$
is so. Hence $\underline{f} :Y^{ss}(L)/\!/G \rightarrow B$ is
flat.
\endproof

\subsection{ Toric setting.}
\begin{toric}
\label{toricsetting}
 Consider again our flat family $f :Y
\rightarrow B$. Let $0 \in B$ be a distinguished point such that
$Y_0$ is a toric variety with the maximal torus ${\mathbb T}$.
Assume in addition that $G$ is a subtorus of ${\mathbb T}$ whose
action on $Y_0$ is induced from that of ${\mathbb T}$.
\end{toric}

Under this assumption, $Y_0^{ss}(L_t)/\!/G$ is a toric variety
with the maximal torus ${\mathbb T}/G$. Hence we have

\begin{cor}
\label{toricgit} Under the toric setting, $\underline{f}
:Y^{ss}(L)/\!/G \rightarrow B$ provides a flat toric degeneration
for generic fibers of $\underline{f}$.
\end{cor}

Apply the above to  toric degenerations of complete and partial
flag varieties, we re-obtain a result of Foth and Hu \cite{FH}
(see also \cite{HMSV05} for the case of $\Gr(2, \C^n)$).

\begin{cor}
\label{toricflag} GIT quotients of complete and partial flag
varieties by the maximal torus degenerate flatly to projective
toric varieties.
\end{cor}

We remark here that the degenerations in Theorem \ref{flatgit},
Corollaries \ref{toricgit} and \ref{toricflag} all depend on the
flat families that we start with.

\section{Flat family of Chow quotients}
\label{chowquotientfamily}

\subsection{ Family of general Chow quotients.} In this section,
we assume in addition that the base $B$ is projective. Hence the
total variety $Y$ is also projective.

As in \cite{Hu2003b}, we will call the Chow cycles in
$Y/\!/^{ch}G$ Chow fibers. By the assumptions of the action, every
Chow fiber is vertical, meaning that it is contained in the fiber
$Y_t$ for some $t \in B$. This leads to a map
$$f^{ch}: Y/\!/^{ch}G \rightarrow B.$$
To see that this is a projective morphism. We may treat $Y$ as a
scheme over $B$ acted upon by the trivial group scheme $G/B=G
\times B$. Use the Chow scheme $\hbox{Ch}(Y/B)$ (\cite{Ko}), we
can embed the space the closures of generic $(G/B)-$ orbits into
$\hbox{Ch}(Y/B)$ and take closure of the image to get the Chow
quotient $(Y/B)/\!/^{ch}(G/B)$. As variety, this is the same as
$Y/\!/^{ch}G$. Hence $f^{ch}: Y/\!/^{ch}G \rightarrow B$ is a
projective morphism. It is trivial to see that the fiber of
$f^{ch}: Y/\!/^{ch}G \rightarrow B$ over $t \in B$ is exactly
$Y_t/\!/^{ch}G$. Hence we have

\begin{thm}
When $B$ is regular and of dimension 1, $f^{ch}$ is a flat family
of Chow quotients. \end{thm}

 We suspect this is always true as long as $f$
is flat. But we do not need this generality for applications.

\begin{thm}
\label{toricChow} Let $f :Y \rightarrow B$ be as in Toric Setting
\ref{toricsetting}. Assume  that the base $B$ is regular and of
dimension 1. Then $f^{ch}$ is a flat family  of the Chow quotients
$Y_t/\!/^{ch}G$ ($t \in B$) with the central fiber $Y_0/\!/^{ch}G$
a projective toric variety.
\end{thm}
\proof It suffices to show that $Y_0/\!/^{ch}G$ is a projective
toric variety, but this follows readily from \cite{KSZ}.
\endproof

\subsection{ Toric settings.}

Now consider any flat toric degeneration of a (partial) flag
variety, e.g, the Grassmannian $\Gr(k, \C^n)$,  over the affine
space $\C$, we can easily extend this family into a flat family
over ${\mathbb P}^1$. Then apply the above discussions to the Chow
quotients of this family, we obtain

\begin{thm}
\label{toricMonbar} For every flat degeneration of a (partial)
flag variety %(e.g, the Grassmannian $\Gr(k, \C^n)$)
 over $\proj^1$
as in Toric Setting \ref{toricsetting} such that the group $G$ is
the maximal torus,  there is an induced flat family over ${\mathbb
P}^1$ such that generic fibers are isomorphic to the Chow quotient
of the flag variety and the central fiber over $0$ is a projective
toric variety.
\end{thm}

Hence, in particular, by applying \cite{LG} and the above, we have

\begin{cor}
The Chow quotients of Grassmannians studied by Kapranov
(\cite{Ka}) and later by Lafforgue (\cite{Lafforgue},
\cite{Lafforgue99}) and  again by Hacking, Keel and Tevelev
(\cite{Hacking03}, \cite{KT}, \cite{HKT}) admit flat toric
degenerations. The same holds for the Chow quotients of all flag
varieties.
\end{cor}

Isolating the special case of $\Gr(2, \C^n)$, we then have the
following important corollary:

\begin{cor}
\label{toricMonbar} For every flat degeneration of $\Gr(2, \C^n)$
over $\proj^1$ as in Toric Setting \ref{toricsetting} such that
the group $G$ is the maximal torus of $\PGL_n$, there is an
induced flat family over ${\mathbb P}^1$ such that generic fibers
are isomorphic to $\overline{M}_{0,n}$ and the central fiber over
$0$ is a projective toric variety $\overline{N}_{0,n}$.
\end{cor}
\proof All we need is to apply Kapranov's theorem that
$\overline{M}_{0,n}$ is isomorphic to the Chow quotient of $\Gr(2,
\C^n)$ by the maximal torus.
\endproof

\begin{rem}
$\overline{N}_{0,n}$ depends on the family that we start with. In
this paper, by abusing notation slightly, we will always use the
same notation $\overline{N}_{0,n}$ to denote any toric
degeneration of $\overline{M}_{0,n}$ that is being considered.
\end{rem}

\section{Toric Chow Quotients of Flag Varieties and Toric  $\overline{M}_{0,n}$}

The results so far obviously apply to: quotients of flag
vareities, quotients of Schubert varieties, quotients of the
products of Grassmannians ( see \cite{Hu2003a}), and so on. But we
will devote the rest of the paper solely to flag varieties and
 $\overline{N}_{0,n}$.
%We will give two descriptions of toric varieties
%$\overline{N}_{0,n}$, both are identified as the Chow quotients of
%some projective toric varieties.

\subsection{ Toric Chow quotients of flag varieties.}
Let $G$ be a connected complex semisimple group\footnote{We note
here that this group $G$ is different than the one used in all the
previous sections. For example, the group $G$ in Toric Setting
\ref{toricsetting} will be replaced instead by a Cartan subgroup
$H$ in this section. We apologize for any possible confusion this
may cause.}, $B$ a Borel subgroup, $U$ its unipotent radical, and
$H$ a Cartan subgroup such that $B=HU$. Let also $\Phi=\Phi(G, H)$
be the system of roots, $\Phi^+=\Phi^+(B, H)$ the subset of
positive roots, and $\{ \alpha_1, ...,\alpha_r\}$ the basis of
simple roots, where $r$ is the rank of $G$. Let $\Lambda$ be the
weight lattice of $G$ and $\Lambda^+$ the subset of dominant
weights. For $\lambda\in\Lambda^+$ we denote by $V(\lambda)$ the
irreducible $G$-module with highest weight $\lambda$. Let
$P_\lambda\supset B$ be the parabolic subgroup of $G$ which
stabilizes a highest weight vector in $V(\lambda)$. Also denote by
$L_\lambda = G \times_{P_\lambda} \C$ the $G$-linearized line
bundle on $X_\lambda:=G/P_\lambda$ corresponding to the character
$\lambda$ extended to $P_\lambda$.

Let $W$ be the Weyl group and $w_0\in W$ the longest element of
length $\ell$. Choose a reduced decomposition
$$
\underline{w_0}=s_{i_1}s_{i_2}\cdots s_{i_\ell}
$$
into a product of simple reflections. The space $A:=\C [G]^U$ of
regular, right $U$-invariant  functions on $G$, has a so-called
 canonical basis $(b_{\lambda, \phi})$, where each
$b_{\lambda, \phi}$ is an eigenvector for both left and right
$H$-action. For the right $H$-action it has weight $\lambda$ and
for the left $H$-action the weight is given by
$$-\lambda+t_1\alpha_{i_1}+\cdots +t_\ell\alpha_{i_\ell}$$
where $b_{\lambda, \phi}$ is parameterized by $$(\lambda, t_1,
..., t_\ell)\in \Lambda^+\times\N^\ell.$$ (The parameterizations
of) $b_{\lambda, \phi}$ generate
 a rational
convex polyhedral cone $${\mathcal C}_{\underline{w_0}}\subset
\Lambda_\R\times \R^\ell$$ (see, e.g., \cite{BZ}). Let $p_1$ and
$p_2$ be the projection of $\Lambda_\R\times \R^\ell$ to the first
and second factor, respectively.  Then for any fixed
$\lambda\in\Lambda^+$,   $Q(\lambda)= p_1^{-1} (\lambda) \cap
{\mathcal C}_{\underline{w_0}}$ is the so-called {\it string
polytope} of $\lambda$.  $Q(\lambda)$ may be identified with its
image in $\R^\ell$ via the second projection $p_2: \Lambda_\R
\times \R^\ell \to \R^\ell$. Let $$ \pi_\lambda:\ \R^\ell\to
\Lambda_\R$$ be defined by  $$(t_1, ..., t_\ell)\mapsto
-\lambda+t_1\alpha_{i_1}+\cdots +t_\ell\alpha_{i_\ell}.$$ Then
$\pi_\lambda$ sends the string polytope $Q(\lambda)$ onto the
convex hull of the Weyl group orbit of the dual weight
$\lambda^*=-w_0 \cdot \lambda$:
$$
\pi_\lambda(Q(\lambda))={\rm Conv}(W \cdot \lambda^*)=-{\rm
Conv}(W \cdot \lambda):=\Delta(\lambda)
$$

Alexeev and Brion (\cite{AlexBrion}) and also Caldero
(\cite{Caldero})  have constructed a flat deformation of the
polarized flag variety $(X_\lambda, L_\lambda)$ to a polarized
toric variety $(X_{\lambda; 0}, L_{\lambda; 0})$ such that the
moment polytope of $(X_{\lambda; 0}, L_{\lambda; 0})$  is  the
string polytope $Q(\lambda)$.

Let $\Phi_\lambda: X_{\lambda; 0} \rightarrow Q(\lambda) \subset
\R^\ell$ be the moment map for ${\mathbb T}$-action, where
${\mathbb T}$ is the compact part of $(\C^*)^\ell$. Then the
composition $\pi_\lambda \circ \Phi_\lambda: X_{\lambda; 0} \to
\Delta_\lambda \subset \Lambda_\R$ is the moment map for the
maximal compact torus $T \subset H$.

Then apply the results from \S \ref{chowquotientfamily}, we have

\begin{thm}
%Let $G = \SL_n$. Choose $\lambda$ such that $X_\lambda = G/P_\lambda$
%is the Grassmannian $\Gr(2,n)$. Then
The Chow quotient of $X_\lambda = G/P_\lambda$ by the maximal
torus $H$ degenerates flatly to the Chow quotient of
 the toric variety $X_{\lambda; 0}$ by the subtorus $H$. The corresponding
 fan of the toric Chow quotient is the common refinement of the normal fans
 of the polytopes $\pi_\lambda^{-1}(\mu) \cap Q(\lambda)$ for all
 $\mu \in \Delta(\lambda)$. In addition, with respect to some suitable ample
 line bundle, this toric Chow quotient corresponds to the fiber
 polytope $\Sigma(Q(\lambda), \Delta_\lambda)$ for the projection
 $\pi_\lambda: Q(\lambda) \rightarrow \Delta_\lambda$.
 \end{thm}
 \proof
The flat family follows from Corollary \ref{toricMonbar}. The rest
is standard for Chow quotients of toric varieties and fiber
polytopes (see \cite{KSZ} and \cite{BS}).
 \endproof

\subsection{The toric variety $\overline{N}_{0,n}$.}
%Next, we offer a more concrete $\overline{N}_{0,n}$.

Let $\Delta_2^n$ be the second hypersymplex defined by
$$\Delta_2^n =\{(r_1, \ldots, r_n) | 0 \le r_i \le 1, \sum_i r_i
=2\}.$$ Any point $(r_1, \ldots, r_n) \in \Delta_2^n$ can be
realized as side lengthes of a polygon in $\R^3$. Let $d_2,
\ldots, d_{n-2}$ be the diagonals emanating from the start of the
first edge to the ends of second to $(n-2)$-th edges. Then $(r_1,
\ldots, r_n, d_2, \ldots, d_{n-2})$ satisfies all possible
triangle inequalities $$d_i - d_{i+1} \le r_{i+1},$$
$$d_i - d_{i+1} \ge -r_{i+1},$$
$$d_i + d_{i+1} \ge r_{i+1}$$
for $1 \le i \le n-2$. Here $d_1=r_1$ and $d_{n-1}=r_n$.  These
inequalities define a polytope $\Gamma^n_2$ in $\R^{2n-3}= \R^n
\times R^{n-3}$, called a Gelfand-Tsetlin polytope (see, e.g.,
\cite{HMSV05} for details). Let $$\varphi: \Gamma^n_2 \rightarrow
\Delta_2^n$$ be the polytopal projection by forgetting the
coordinates $d_2, \ldots, d_{n-2}$. Let $\Sigma(\Gamma^n_2,
\Delta_2^n)$ be the fiber polytope of the projection $\varphi:
\Gamma^n_2 \rightarrow \Delta_2^n$ (\cite{BS}).

\begin{thm}
 $\overline{M}_{0,n}$ admits a
flat degeneration to the toric variety $\overline{N}_{0,n}$
defined by the fiber polytope $\Sigma(\Gamma^n_2, \Delta_2^n)$.
The fan of $\overline{N}_{0,n}$ is the common refinement of the
normal fans of all the polytopes $\varphi^{-1}((r_1, \ldots,
r_n)), (r_1, \ldots, r_n) \in \Delta_2^n$.
\end{thm}
\proof Take a flat family that degenerates $\Gr(2,\C^n)$ to the
toric variety defined by $\Gamma^n_2$ (see, e.g., \cite{Caldero},
\cite{LG} for the construction), then apply the results from \S
\ref{chowquotientfamily} and \cite{KSZ}.
\endproof

We point it out again that $\overline{N}_{0,n}$ is sensitive to
many choices involved. The one in the theorem is just a particular
example of {\it many} toric degenerations of $\overline{M}_{0,n}$

\begin{rem}
In \cite{KT05}, Keel and  Tevelev wrote down some equations for
$\overline{M}_{0,n}$. Ideally, it would be nice if one can produce
some binomials out of their equations, giving rise to certain
toric degenerations of $\overline{M}_{0,n}$ with interesting
binomial equations.
\end{rem}

\subsection{Moduli interpretations of $\overline{N}_{0,n}$}

In some symplectic way, an interpretation can be done in terms of
stable polygons (\cite{Hu99}). Topologically, this toric
$\overline{M}_{0,n}$ can be obtained by collapsing the moduli
spaces of stable polygons so that the bending flows extend
everywhere (\cite{Huinprepare}). This is similar to Kamiyama and
Yoshida's construction for the moduli space of ordinary polygons
(\cite{KY}). But, we would be happier to have an algebro-geometric
moduli interpretation.

Further details of these toric  $\overline{N}_{0,n}$, including
their moduli as well as polygonal interpretations,  will appear in
forthcoming papers (\cite{Huinprepare}).

\medskip
{\sl Acknowledgements.} Many works have inspired this paper, in
particular all those interesting papers on toric degenerations of
flag varieties, Kapranov's paper that realizes
$\overline{M}_{0,n}$ as Chow quotient, and
Kapranov-Sturmfels-Zelevinsky's paper on (Chow) quotients of toric
varieties. The paper of Alexeev and Brion (\cite{AlexBrion}) on
toric deformation has been particularly helpful. I was convinced
of toric degenerations of $\overline{M}_{0,n}$ for quite a while,
I thank John Millson whose visit  in the Spring of 2004
%(during the Arizona Conference on the Geometry and Topology of Quotients)
convinced me that it should be written up.

%\bibliographystyle{amsplain}
%\makeatletter \renewcommand{\@biblabel}[1]{\hfill#1.}\makeatother
%\newcommand{\bysame}{\leavevmode\hbox to3em{\hrulefill}\,}

\end{document}